\documentclass[tbtags,reqno]{amsart}

\usepackage{young}
  
\makeatletter
\renewcommand{\subsubsection}{\@startsection
{subsubsection}
{3}
{0mm}
{\baselineskip}
{-0.5\baselineskip}
{\normalfont\normalsize\bfseries}}
\makeatother

\newtheorem{theorem}{Theorem}
\newtheorem{lemma}[theorem]{Lemma}

\theoremstyle{definition}

\theoremstyle{remark}

\newtheorem*{acknow}{Acknowledgments}

\begin{document}

 \title[A short proof of the integrality of the ($q,t$)-Kostka coefficients]{A short proof of the integrality of the Macdonald ($q,t$)-Kostka coefficients}
\author{Luc~Lapointe and Luc~Vinet}
\address
{Centre de Recherches Math\'ematiques, \\Universit\'e de
Montr\'eal,  C.P.~6128,  succursale~Centre-ville, \\ Montr\'eal,
Qu\'ebec, Canada, H3C 3J7}
 \date{}

\begin{abstract}
The Macdonald polynomials can be obtained by acting on the constant 1 with creation operators.  Three different expressions for these operators are derived, one from the other, in a rather succint way.  When the last of these expressions is used, the formalism is seen to imply straightforwardly the integrality of the ($q,t$)-Kostka coefficients, that is of the expansion coefficients for the Macdonald functions in terms of Schur functions.
\end{abstract}
\maketitle

\section{Introduction and background}

Let $\Lambda_N$ denote the ring of symmetric functions in the
variables $x_1,x_2,\dots,x_N$ and denote by $\mathbb Q(q,t)$ the
field of rational functions of the parameters $q$ and $t$ with rational
coefficients.  The Macdonald polynomials \cite{1} $J_{\lambda}(x;q,t)$ are symmetric polynomials labelled by
partitions $\lambda=(\lambda_1,\lambda_2,\ldots)$, where $\lambda_1
\geq \lambda_2 \geq \dots$.  i.e.  sequences of non-negative integers 
in decreasing order.  These polynomials form a basis for $\Lambda_N
\otimes \mathbb Q(q,t)$ and can be characterized as the joint
eigenfunctions of the commuting operators $\{M_N^k,
k=0,\dots,N\}$ defined as follows:
\begin{equation}
 M_N^k=\sum_I t^{(N-k)k+k(k-1)/2}\tilde{A}_I(x;t)\prod_{i\in I}
T_{q,x_i} 
\end{equation}
with $M_N^0=1$.  Here, the sum goes over all $k$-element subsets  $I$ of
$\{1,\dots,N\}$,
\begin{equation}
\tilde{A}_I (x;t)=\prod_{\begin{subarray} {c} i\in I\\ j\in\{1,\ldots,N\}
\backslash I \end{subarray}} \frac{x_i-t^{-1}x_j}{  x_i-x_j} 
\end{equation}
and $T_{q,x_i}$ stands for the $q$-shifted operator in the variables $x_i$ ($T_{q,x_i} f(x_1,\dots,x_i,\dots)$ $=$ $f(x_1,\dots,q x_i,\dots)$).  This notation will be used throughout the
paper.  The eigenvalue equations that the Macdonald polynomials satisfy are conveniently written in terms of the generating function
\begin{equation}
M_N (X;q,t) = \sum_{k=0}^N M_N^k X^k,
\end{equation}
where $X$ is an arbitrary parameter. With $J$, a set of cardinality $|J|=j$, we  shall also  use the notation $M_J (X;q,t)$ to designate $M_{j}(X;q,t)$ in the variables $x_i$, $i \in J$.  If $\ell(\lambda)$ denotes the number of non-zero parts of $\lambda$, for  $\ell(\lambda) \leq N$ we have
\begin{equation}
M_N(X;q,t) J_{\lambda}(x;q,t) = a_{\lambda}(X;q,t) J_{\lambda}(x;q,t),
\end{equation}
with 
\begin{equation}
a_{\lambda}(X;q,t) = \prod_{i=1}^N (1+X q^{\lambda_i} t^{N-i}).
\end{equation}
It is customary to denote by $a(s)$ and $\ell(s)$ the number of squares in the diagram of $\lambda$ that are respectively to the south and east of the square $s$ in the diagram of the partition $\lambda$.  A complete characterization of the Macdonald polynomials can now be given by supplementing (4) with the following decomposition of $J_{\lambda}(x;q,t)$ on the monomial symmetric functions $m_{\mu}$:
\begin{equation}
 J_\lambda = \sum_{\mu \le \lambda} v_{\lambda\mu}(q,t) m_\mu, 
\end{equation}
with
\begin{equation}
 v_{\lambda\lambda}(q,t) \equiv
c_{\lambda}(q,t) = \prod_{s \in \lambda} (1-q^{a(s)} t^{\ell (s)+1}).
\end{equation}
The sum in (6) is over partitions that are smaller than (or equal to) $\lambda$ in  the dominance order.  Recall that $m_{\mu} = \sum_{\text{distinct permutations}} x_1^{\mu_1} x_2^{\mu_2} \dots $.

We are looking for creation operators $B_k$, $k=1,\dots,N$, such that, for any partition $\lambda$ with $\ell(\lambda) \leq k$,
\begin{equation}
B_k J_{\lambda}(x) = J_{\lambda+(1^k)}(x).
\end{equation}
Given such creation operators, it is clear that the Macdonald polynomials associated to any partition can be constructed by acting repeatedly with these $B_k$ on $J_{0} (x)=1$.  Indeed, the following Rodrigues formula for $J_{\lambda}(x)$ is an immediate consequence of (8):
\begin{equation} 
J_\lambda(x)
	= B_N^{\lambda_N} B_{N-1}^{\lambda_{N-1} -\lambda_N} \dots B_1^{\lambda_1 - \lambda_2} \cdot 1 .
\end{equation}

Three expressions $B_k^{(i)}$, $i=1,2,3$, have been obtained for the creation operators.  ( In fact, we also know a fourth that will not be used here; see eq. (43)-(44) in \cite{2}).  When $B_k^{(3)}$ is used, formula (9) is readily seen to imply \cite{3,4} the integrality of the $(q,t)$-Kostka coefficients $K_{\lambda \mu}(q,t)$.  The purpose of this note is to provide a simple derivation of these expressions and therefore,  a short proof  of the fact that the 
$K_{\lambda \mu}(q,t)$ are polynomials in $q,t$ with integral coefficients.

We need more notation \cite{5} to write down the expressions $B_k^{(1)}$, $B_k^{(2)}$ and $B_k^{(3)}$.  For $n$ an integer, the $q$-shifted factorials are defined by
\begin{equation}
 (a;q)_n= (1-a)(1-qa) \dots (1-q^{n-1}a); \quad  (a;q)_0 \equiv 1.
\end{equation}
They satisfy various identities among which:
\begin{equation}
\begin{split}
(a;q)_{m+n} & = (a;q)_m (a q^m;q)_n\\
(a;q)_{n-k} & = \frac{(a;q)_n }{(q^{1-n}/a;q)_k} \Bigl( -\frac{q}{a} \Bigr)^k q^{k(k-1)/2-n k}.
\end{split}
\end{equation}
The $_1\phi_1$ basic hypergeometric series is
\begin{equation}
_1\phi_1 (a;b;q,x) = \sum_{n=0}^{\infty} (-1)^n q^{n(n-1)/2} \frac{(a;q)_{n}}{(q;q)_{n}(b;q)_{n}} x^n;
\end{equation}
it can be summed when $x=b/a$ and, in particular, we have
\begin{equation}
_1\phi_1 (t^{-n};t^{-n+1}q^{-1};t,t/q)=\frac{1}{(t^{-n+1}q^{-1};t)_n}.
\end{equation}

At some point in the proofs, we shall also need another summation formula, this one a special case of a result due to Garsia and Tesler (\cite{6}, Proposition~3.1): let $J \subseteq \{1,\dots,N\}, J^c=\{1,\dots,N\} \backslash J$ and $x_J= \prod_{j \in J} x_j$, for all $k$ and $m$ we have
\begin{equation}
\sum_{|J|=k} x_J \sum_{\begin{subarray}{c} J' \subseteq J^c \\ |J'|=m \end{subarray}}  \tilde A_{J \cup J'}=t^{-m(N-k-m)} 
\begin{bmatrix}
N-k \\
m
\end{bmatrix}_t
 \sum_{|J|=k} x_J,
\end{equation}
where $\tilde A_{J \cup J'}$ is a defined in (2) and where the second sum on the l.h.s. of (14) is over all subsets $J'$ of $J^c$ with cardinality $|J'|=m$.  The $q$-binomial coefficient is 
\begin{equation}
\begin{bmatrix}
n\\
k
\end{bmatrix}_q
=\frac{(q;q)_{n}}{(q;q)_{k}(q;q)_{n-k}}.
\end{equation}

We can now give the three expressions $B_{k}^{(i)}$ $i=1,2,3$, $k=1,\dots,N$ of the creation operators that we will derive in the next section.

\noindent $\bullet$ Expression 1
\begin{equation}
B_k^{(1)} = \frac{1}{(q^{-1};t^{-1})_{N-k}} M_N(-t^{k+1-N}q^{-1};q,t) e_k.
\end{equation}
Here $e_k$ stands for the $k^{th}$ elementary function
\begin{equation}
e_k(x)=\sum_{|I|=k} x_I=m_{(1^k)}.
\end{equation}
 $\bullet$ Expression 2
\begin{equation}
B_k^{(2)}=\sum_{|I|=k}x_I  \sum_{m=0}^{N-k} \sum_{\begin{subarray}{c} I' \subseteq I^c \\ |I'|=m \end{subarray}} \frac{q^{-m}}{(t^{k+1-N}q^{-1};t)_m} \tilde A_{I \cup I'} M_{I \cup I'}(-t^{1-m};q,t).
\end{equation}
 $\bullet$ Expression 3
\begin{equation}
B_k^{(3)} = \sum_{|I|=k} x_I \tilde A_I M_{I}(-t;q,t).
\end{equation}

Our proofs of these formulas will proceed like this.  The first expression for the creation operators will be seen to follow from the Pieri formula which gives the action of the elementary functions $e_k$ on the monic Macdonald polynomials $P_{\lambda}=1/c_{\lambda} (q,t) J_{\lambda}$.  This formula reads \cite{1}
\begin{equation}
e_k P_{\lambda}= \sum_{\mu} \Psi_{\mu/\lambda} P_{\mu},
\end{equation}
where the sum is over all partitions $\mu$ containing $\lambda$ such that the set-theoretic difference $\mu-\lambda$ is $k$-dimensional with the property that $\mu_i -\lambda_i \leq 1$, $ \forall i \geq 1$.  If  $C_{\mu/\lambda}$  and  $R_{\mu/\lambda}$ respectively  denote the union of the columns and of the rows  that intersect $\mu - \lambda$, the coefficients $\Psi_{\mu/\lambda}$ are given by
\begin{equation} 
\Psi_{\mu/\lambda}= \prod_{\begin{subarray}{c} s \in C_{\mu/\lambda} \\ s 
\not \in R_{\mu/\lambda} \end{subarray}} \frac{b_{\mu}(s)}{b_{\lambda}(s)}
\end{equation}
where
\begin{equation}
b_{\lambda}(s) = {\cases \frac{1-q^{a(s)} t^{\ell(s)+1}}{1-q^{a(s)+1} t^{\ell(s)}} \qquad & {\text{if}}~s \in \lambda \\
1    & {\text{otherwise}} \endcases}.
\end{equation}
Second, we shall show that $B_k^{(2)}=B_k^{(1)}$ with the help of the summation formula (13) and (14).  The third expression for the creation operators will then be easily arrived at using the second one by observing that  $B_k^{(2)}$ reduces to  $B_k^{(3)}$ when acting on Macdonald polynomials; in other words, we shall prove that  $B_k^{(2)} J_{\lambda}=B_k^{(3)} J_{\lambda}= J_{\lambda+(1^k)}$ for $\ell(\lambda) \leq k$.  Creation operators of type 3 were first obtained in \cite{7} in the limit case of the Jack polynomials  $q=t^{\alpha},t \to 1$.   A realization of these operators in terms of Dunkl operators was given and the Rodrigues formula they entail implied the integrality of the expansion coefficients of the Jack polynomials over the symmetric monomial basis \cite{8}.

In the case of the Macdonald polynomials, Expression 1 was first derived in \cite{2};  Expression 3 is also given there in the form of a conjecture.  Meanwhike Kirillov and Noumi provided a realization of $B_k^{(3)}$ in terms of Dunkl-Cherednik operators  \cite{3,4} and also gave two proofs of the fact that they are indeed creation operators for the Macdonald polynomials.  Expression 2 seems new; it plays an important role as the essential intermediate step in the concise derivation of $B_k^{(3)}$ from $B_k^{(1)}$ that we present next.

In a separate paper \cite{9}, we present another  derivation of  $B_k^{(3)}$.  The starting point is again Expression 1, i.e. $B_k^{(1)}$.  The approach followed in \cite{9} differs from the one adopted in the present one; it uses in an essential way the specifics of the realization in terms of Dunkl-Cherednik operators and properties of affine Hecke algebras.  The derivation is more involved and lengthy than the one given here.  It offers however the advantage of being more constructive allowing one to basically arrive at Expression
2 from Expression 1.  Similar developments in the case of the Jack
polynomials deserve special attention and are the object of a forthcoming
publication \cite{10}.

\section{Proofs}

We now proceed along the lines indicated in the introduction to show
in turn that $B^{(1)}_k$, $B^{(2)}_k$ and $B^{(3)}_k$ act as creation
operators for the Macdonald polynomials.

\begin{theorem}
  For any partition $\lambda$ with $l(\lambda) \leq k$,
the operators $B^{(1)}_k$ act as follows on the Macdonald polynomials
$J_{\lambda}(x;q,t)$:
\begin{equation}
B^{(1)}_k J_{\lambda}(x) = J_{\lambda+(1^k)}(x).
\end{equation}
\end{theorem}
Proof.  The following lemma is an immediate consequence of the Pieri formula.
\begin{lemma} For $\lambda$ a partition such that $\ell(\lambda) \leq k$, the action of $e_k$ on $P_{\lambda}$ is given by 
\begin{equation}
e_k P_{\lambda} =  P_{\lambda+(1^k)} +
\sum_{\mu \neq \lambda+(1^k)} \Psi_{\mu/\lambda} P_{\mu},
\end{equation} 
where all the $\mu$'s in the sum are such that $\mu_{k+1}=1$.
\end{lemma}
Indeed, the only way to construct a $\mu$ with $\mu_{k+1} \neq 1$ is to add a 1 in each of the first $k$ entries of $\lambda$.
From Lemma~2 and (4) and (5), we have
\begin{equation}
 M_N(-t^{k+1-N}q^{-1};q,t) e_k  P_{\lambda} = \prod_{i=1}^k (1-t^{k+1-i}q^{\lambda_i}) (q^{-1};t^{-1})_{N-k} P_{\lambda+(1^k)}
\end{equation}
since the eigenvalues
\begin{equation}
a_{\mu}(-t^{k+1-N}q^{-1};q,t) = \prod_{i=1}^N (1-t^{k+1-i}q^{\mu_i-1}),
\end{equation}
of  $ M_N(-t^{k+1-N}q^{-1};q,t)$ on the $P_{\mu}$'s in (24) vanish if $\mu_{k+1}=1$.

From the definition given in (7), it is easy to check that
\begin{equation}
\frac{c_{\lambda+(1^k)}}{c_{\lambda}}=\prod_{i=1}^k (1-t^{k+1-i}q^{\lambda_i}).
\end{equation}
Using this result and passing from $P_{\lambda}$ to $J_{\lambda}$ we see that 
\begin{equation}
B_k^{(1)}J_{\lambda}= \frac{1}{(q^{-1};t^{-1})_{N-k}} M_N(-t^{k+1-N}q^{-1};q,t) e_k  J_{\lambda}= J_{\lambda+(1^k)}
\end{equation}
when $\ell(\lambda) \leq k$.  This proves Theorem~1.
\begin{theorem}
\begin{equation}
B_k^{(2)}=B_k^{(1)},
\end{equation}
hence $B_k^{(2)}$ is also such that $B_k^{(2)} J_{\lambda}=J_{\lambda+(1^k)}$ for $\ell(\lambda) \leq k$.
\end{theorem}
Proof.  We first prove the following lemma.
\begin{lemma} Let $J$ and $\bar J$ be complementary subsets of $\{1,\dots,N\} \backslash \{1,\dots,\ell\}$, $\ell \leq N$.  For all $n,k$ such that $N-\ell \geq k-n \geq 0$,
\begin{equation}
\begin{split}
\sum_{m'=0}^{N-k-\ell +n} \frac{q^{-m'}}{(t^{k-N+1}q^{-1};t)_{m'+\ell-n}} 
& \sum_{|J|=k-n} x_J \\ 
\times \sum_{\begin{subarray}{c} J' \subseteq \bar J \\ |J'|=m' \end{subarray}} & \tilde A_{J \cup J'}^{\{1,\dots,N\} \backslash \{1,\dots,\ell\}} = 
 \frac{1}{(q^{-1};t^{-1})_{N-k}} \sum_{|J|=k-n} x_J,
\end{split}
\end{equation}
where 
\begin{equation}
\tilde A_{I}^J =\prod_{\begin{subarray}{c} i \in I \\ j \in J \backslash I  \end{subarray}} \frac{x_i - t^{-1}x_j}{x_i -x_j}.
\end{equation}
\end{lemma}
To derive this relation, one first uses formula (14) in the l.h.s. of (30).  After simplifying the factors $\sum_{|J|=k-n}x_J$, one thus finds that (30) amounts to
\begin{equation}
\sum_{m'=0}^{N-k-\ell +n} \frac{q^{-m'} t^{-m'(N-k-\ell +n -m')}
}{(t^{k-N+1}q^{-1};t)_{m'+\ell-n}}\begin{bmatrix}
N-k-\ell+n\\
m'
\end{bmatrix}_t 
= \frac{1}{(q^{-1};t^{-1})_{N-k}} .
\end{equation}
With the help of the identities (11), one easily shows that
\begin{equation}
{\text{l.h.s. of (32)}}=\frac{1}{(t^{k-N+1}q^{-1};t)_{\ell-n}} {_1\phi_1} (t^{k-N+\ell-n};t^{k-N+\ell-n+1}q^{-1};t,t/q).
\end{equation}
At this point, we recall the $_1\phi_1$ sum in (13) to find that
\begin{equation}
{\text{l.h.s. of (32)}}=\frac{1}{(t^{k-N+1}q^{-1};t)_{\ell-n}} \frac{1}{(t^{k-N+\ell-n+1}q^{-1};t)_{N-k-\ell+n}}
\end{equation}
and to see using (11) that (30) is indeed verified.

We now return to the proof of (29), that is of the identity $B_k^{(2)}=B_k^{(1)}$.  Let $\ell \leq N$.  Since $B_k^{(1)}$ and $B_k^{(2)}$ are symmetric, it will suffice to show that the coefficient to the left of the operators $T_{q,x_{1}} \dots T_{q,x_{\ell}}$ in $B_k^{(1)}$ and $B_k^{(2)}$  are identical.

To that end, let 
\begin{equation}
\begin{split}
& I=L \cup J, \quad \quad I'=\bar L \cup J',\\
& L \subseteq \{1,\dots, \ell\}, \quad \quad \bar L = \{1,\dots,\ell\} \backslash L,\\
& J,J' \subseteq \{1,\dots,N\} \backslash \{1, \dots,\ell \}= J \cup \bar J,\\
& J \cap \bar J= \phi, J' \subseteq \bar J,
\end{split}
\end{equation}
and define
\begin{equation}
[\ell,k]= {\cases \ell \qquad & {\text{if}}~\ell \leq k \\
k    &  {\text{if}}~\ell > k \endcases}.
\end{equation}
Collecting the factors to the left of $T_{q,x_{1}}\dots T_{q,x_{\ell}}$ in $B_k^{(1)}$ and $B_k^{(2)}$ ( see (16) and (18)) gives respectively:
\begin{equation}
\begin{split}
B_k^{(1)}|_{T_{q,x_{1}}\dots T_{q,x_{\ell}}}= 
 \frac{1}{(q^{-1};t^{-1})_{N-k}} & \bigl( -t^{k+1-N}q^{-1}\bigr)^{\ell}\\
\times \tilde A_{\{1,\dots,\ell\}} & t^{\ell(\ell-1)/2} t^{\ell(N-\ell)} \sum_{n=0}^{[\ell,k]}  \sum_{|L|=n} q^n x_L \sum_{|J|=k-n} x_J
\end{split}
\end{equation}
and
\begin{equation}
\begin{split}
 B_k^{(2)}|_{T_{q,x_{1}}\dots T_{q,x_{\ell}}}&= 
 \sum_{n=0}^{[\ell,k]}   \sum_{|L|=n}  x_L \sum_{|J|=k-n} x_J \sum_{m=\ell-n}^{N-k}  \sum_{|\bar L \cup J'|=m}\\
 \times &  \frac{q^{-m}}{(t^{k+1-N}q^{-1};t)_{m}} \bigl( -t^{1-m}\bigr)^{\ell}  t^{\ell(\ell-1)/2} t^{\ell(m+k-\ell)} \tilde A_{J \cup J'}^{J \cup \bar J} \tilde A_{\{1,\dots,\ell\}},
\end{split}
\end{equation}
where we have used in (38) the fact that $  \tilde A_{I \cup I'} \tilde A_{\{1,\dots,\ell\}}^{I \cup I'}=\tilde A_{J \cup J'}^{\{1,\dots,N\} \backslash \{1, \dots,\ell \}} \tilde A_{\{1,\dots,\ell\}} $.  It is then immediate to see that the equality
\begin{equation}
B_k^{(1)}|_{T_{q,x_{1}}\dots T_{q,x_{\ell}}}=B_k^{(2)}|_{T_{q,x_{1}}\dots T_{q,x_{\ell}}}
\end{equation}
holds, since after trivial simplifications it is seen to amount to
\begin{equation}
\begin{split}
& \sum_{n=0}^{[\ell,k]}   \sum_{|L|=n} q^{n-\ell } x_L \biggl( \frac{1}{(q^{-1};t^{-1})_{N-k}} \sum_{|J|=k-n} x_J \biggr)=\\
& \sum_{n=0}^{[\ell,k]}   \sum_{|L|=n} q^{n-\ell} x_L \biggl( \sum_{m'=0}^{N-k-\ell+n} \frac{q^{-m'}}{(t^{k-N-1}q^{-1};t)_{m'+\ell-n}} \sum_{|J|=k-n} x_J \sum_{|J'|=m'} \tilde A_{J \cup J'}^{J \cup \bar J} \biggr)
\end{split}
\end{equation}
and hence to follow from Lemma~4.

Once Theorem~3 is proved, the main result, Theorem~5, is readily obtained.
\begin{theorem}
For any partition $\lambda$, such that $\ell (\lambda) \leq k$, the actions of $B_k^{(2)}$ and $B_k^{(3)}$ on the Macdonald polynomials $J_{\lambda}(x)$ coincide:
\begin{equation}
B_k^{(3)}J_{\lambda}(x)=  B_k^{(2)} J_{\lambda}(x)= J_{\lambda+(1^k)}(x).
\end{equation}
\end{theorem}
This is shown to be true with the help of the following lemma
\begin{lemma} 
Let  $|I|=k$ and $|I'|=m$, $I' \subseteq I^c$.
\begin{equation}
M_{I \cup I'} (-t^{1-m};q,t) J_{\lambda}(x;q,t)=0,
\end{equation}
if $\ell(\lambda) \leq k $ and  $m > 0$.
\end{lemma}
Proof.  Denote by $x(I)$  the set of variables $\{ x_i, i \in J\}$.  The Macdonald polynomials are known to enjoy the property according to which
\begin{equation}
J_{\lambda} \bigl( x(I),x(I^c) \bigr)  = \sum_{\mu,\nu} \tilde f_{\mu \nu}^{\lambda} J_{\mu} \bigl( x(I) \bigr) J_{\nu} \bigl( x(I^c) \bigr)
\end{equation}
with $ \tilde f_{\mu \nu}^{\lambda}=0$ unless $ \mu \subset \lambda$ and  $ \nu \subset \lambda$ and in particular if $\ell(\mu)$ or $\ell(\nu)$ is greater than $k$. 

  Since $M_{I \cup I'} (-t^{1-m};q,t)$ is a $q$-difference operator depending only of  the variables $x_i$, $i \in I \cup I'$, we see from (43) that 
\begin{equation}
M_{I \cup I'} J_{\lambda}(x)  = \sum_{\mu,\nu} \tilde f_{\mu \nu}^{\lambda} J_{\mu} \bigl( x((I \cup I')^c) \bigr) M_{I \cup I'} J_{\nu} \bigl( x(I \cup I') \bigr).
\end{equation}
The proof of Lemma~6 is then completed by observing from (4)-(5) that
\begin{equation}
M_{I \cup I'} (-t^{1-m};q,t)  J_{\nu} \bigl( x(I \cup I') \bigr) = \prod_{i=1}^{k+m} (1-q^{\nu_i}t^{k+1-i}) J_{\nu} \bigl( x(I \cup I') \bigr)=0,
\end{equation}
whenever  $m>0$, since $\nu_{k+1}=0$.

Theorem~5 is thus an immediate consequence of this lemma since in the expression (18) of $B_k^{(2)}$, all the terms, except the $m=0$ one, are seen to act trivially on the $J_{\lambda}(x)$ when $\ell(\lambda) \leq k$.

\section{Integrality of the $(q,t)-$Kostka coefficients and conclusion}

The $(q,t)-$Kostka $K_{\lambda\mu}(q,t)$ coefficients give the overlaps
between the Macdonald and the Schur polynomials (see \cite{1}):
\begin{equation}
J_{\mu}(x;q,t) = \sum_{\lambda} K_{\lambda \mu} (q,t) S_{\lambda}(x;t).
\end{equation}
As mentioned in the introduction, a proof of Theorem 5 readily
provides a proof of the fact that these $K_{\lambda\mu}(q,t)$'s 
are polynomials in $q$ and $t$ with integer coefficients as was
first conjectured by Macdonald.  How such a conclusion is arrived at 
from Theorem 5 is explained in \cite{3,4}.  The main observation (see also
\cite{8}) is that Theorem 5 implies that the expansion coefficients
$v_{\lambda\mu}(q,t)$ in (6) are polynomials in $q,t$ with integer
coefficients.  One then recalls \cite{1} that the overlaps between the 
monomial symmetric functions and the Schur functions have the form
$p(t)/q(t)$, with $p(t)$ and $q(t)$ polynomials in $t$ with
integer coefficients and such that $q(0)=1$.  Upon using the duality
properties of the function $J_{\lambda}$ and $S_{\lambda}$ one then
straightforwardly concludes that $K_{\lambda,\mu}(q,t) \in 
\mathbb{Z}[q,t]$.

Let us mention in concluding that a other proofs of the
integrality of the $(q,t)$-Kotska coefficients have been given
recently using different approaches by  Garsia and Remmel \cite{11}, Garsia
and Tesler \cite{6}, Knop \cite{12,13} and  Sahi \cite{14}.

\begin{acknow}
We would like to express our thanks to Adriano Garsia  for various comments and suggestions.  

\noindent This work has been supported in part through funds provided by NSERC (Canada) and FCAR (Qu\'ebec).  L.~Lapointe holds a NSERC postgraduate scholarship.
\end{acknow}

\end{document}